\documentclass[preprint,3p,authoryear]{elsarticle}
\journal{Statistics and Probability Letters}
\usepackage{amsthm,amssymb,amsmath}
\newtheorem{theorem}{Theorem}

\newtheorem{corollary}{Corollary}

\newtheorem{lemma}{Lemma}
\theoremstyle{remark}
\newtheorem{remark}{Remark}

\newcommand{\F}{\mathcal{F}}
\renewcommand{\Pr}{\mathbb{P}}
\newcommand{\E}{\mathbb{E}}
\newcommand{\D}{\mathbb{D}[0,\infty)}
\newcommand{\DT}{\mathbb{D}[0,T]}
\newcommand{\DDD}{\mathbb{D}[0,1]}

\begin{document}
\begin{frontmatter}

\title{Functional limit theorems for linear processes in the domain of attraction of stable laws}
\author{Marta Tyran-Kami\'nska}
\ead{mtyran@us.edu.pl}
\address{Institute of Mathematics, University of Silesia, Bankowa 14, 40-007 Katowice,
Poland}

\begin{abstract} We study functional limit theorems for linear type
processes with short memory under the assumption that the innovations are dependent identically distributed
random variables with infinite variance and in the domain of attraction of stable laws.
\end{abstract}
\begin{keyword}
L\'evy stable processes\sep functional limit theorem\sep Skorohod topologies\sep $\rho$-mixing

\MSC 60F05\sep 60F17 
\end{keyword}


\end{frontmatter}
\section{Introduction}

We consider the linear process $\{Z_{j}\colon j\in\mathbb{Z}\}$ defined by
\begin{equation}\label{e:map}
Z_j=\sum_{k=-\infty}^\infty a_k \xi_{j-k},
\end{equation}
where the innovations $\{\xi_j\colon j\in\mathbb{Z}\}$ are identically distributed random variables with
infinite variance and the sequence of constants $\{a_k\colon k\in\mathbb{Z}\}$ is such that $
\sum_{k\in\mathbb{Z}} |a_k|<\infty$. This case is referred to as short memory, or as short range dependence. The
functional central limit theorem (FCLT) for the partial sums of the linear process, properly normalized,  merely
follows from the corresponding FCLT for the innovations being in the domain of attraction of the normal law, see
\cite{peligradutev06b}. Then the limiting process has continuous sample paths and choosing the right topology in
the Skorohod space $\DDD$ is not problematic. However, as shown by~\cite{avramtaqqu92}, the weak convergence of
the partial sums of the linear process with independent innovations (i.i.d.~case) in the domain of attraction of
non-normal laws is impossible in the Skorohod $J_1$ topology on $\DDD$, but the functional limit theorem might
still hold, under additional assumptions, in the weaker Skorohod $M_1$ topology (see \cite{skorohod56}).
\cite{avramtaqqu92} use the standard approach through tightness plus convergence of finite dimensional
distributions. Here, we use approximation techniques and study weak convergence in $\D$, i.e. the space of
functions on $[0,\infty)$ that have finite left-hand limits and are continuous from the right. Given processes
$X_n$, $X$ with sample paths in $\D$, we will denote by $X_n(t)\Longrightarrow X(t)$ the weak convergence in
$\D$ with one of the Skorohod topologies $J_1$ or $M_1$,  and write $\overset{J_1}{\Longrightarrow}$ or
$\overset{M_1}{\Longrightarrow}$, if the indicated topology is used. Note that if the limiting process $X$  has
continuous sample paths then $\Longrightarrow$ in $\D$ with one of the Skorohod topologies is equivalent to weak
convergence in $\D$ with the local uniform topology. For definitions and properties of the topologies we refer
to \cite{jacodshiryaev03} and \cite{whitt02}.

To motivate our approach we first consider the linear process $\{Z_j\colon j\in\mathbb{Z}\}$ as in
\eqref{e:map}, where $\{\xi_j\colon j\in\mathbb{Z}\}$ is a sequence of i.i.d.~random variables. There exist
sequences $b_n>0$ and $c_n$ such that the partial sum processes of the i.i.d. sequence $\{\xi_j\colon
j\in\mathbb{Z}\}$ converge weakly to an $\alpha$-stable L\'evy process $X$ with $0 < \alpha< 2$ \citep[see,
e.g.][Proposition~3.4]{resnick86}
\begin{equation}\label{eq:cj1}
\frac{1}{b_n}\sum_{j=1}^{[nt]}(\xi_j-c_n)\overset{J_1}{\Longrightarrow}X(t)
\end{equation}
if and only if there is convergence in distribution
\[
\frac{1}{b_n}\sum_{j=1}^{n}(\xi_j-c_n)\overset{d}{\longrightarrow}X(1)\quad \text{in}\quad\mathbb{R},
\]
or, equivalently, there exist $p\in[0,1]$ and a slowly varying function $L$, i.e., $L(sx)/L(x)\to 1$ as
$x\to\infty$ for every $s>0$, such that
\begin{equation}\label{e:xitail1}
\lim_{x\to\infty}\dfrac{\Pr(\xi_1>x)}{\Pr(|\xi_1|>x)}=p\quad \mbox{and}\quad \lim_{x\to\infty}\frac{
\Pr(|\xi_1|>x)}{x^{-\alpha}L(x)}=1;
\end{equation}
in that case the sequences $\{b_n,c_n\colon n\in\mathbb{N}\}$ in \eqref{eq:cj1} can be chosen as
\begin{equation}\label{d:bncn}
b_n=\inf\bigl\{x:\Pr(|\xi_1|\le x)\ge 1-\frac1n\bigr\} \quad \text{and}\quad c_n=\E(\xi_1I(|\xi_1|\le b_n)).
\end{equation}
We have $b_n\to\infty$, $nP(|\xi_1| > b_n) \to 1$, and $nb_n^{-\alpha}L(n)\to 1$, as $n\to\infty$, by
\eqref{e:xitail1}. We refer to \cite{feller71} for $\alpha$-stable random variables and their domains of
attraction. If $\alpha=2$, condition \eqref{eq:cj1} holds  if and only if the function $x\mapsto
\E(\xi_1^2I(|\xi_1|\le x))$ is slowly varying; in that case the sequence $b_n$ can be chosen as satisfying
$nb_n^{-2}\E(\xi_1^2I(|\xi_1|\le b_n))\to 1$ as $n\to\infty$, the $c_n$ as in \eqref{d:bncn}, and $X$  is
a~Brownian motion.   In any case, \eqref{eq:cj1} implies that the function $x\mapsto \E(\xi_1^2I(|\xi_1|\le x))$
is regularly varying with index $2-\alpha$, that is, there exists a slowly varying function $\ell$ such that (if
$\alpha<2$ then $\ell(x)=\alpha L(x)/(2-\alpha)$ where $L$ is as in~\eqref{e:xitail1})
\begin{equation}\label{a:daa}
\lim_{x\to\infty}\frac{\E (\xi_1^2I(|\xi_1|\le x))}{x^{2-\alpha}\ell(x)}=1.
\end{equation}
\cite{astrauskas83} and \cite{davisresnick85} show that if the coefficients $\{a_k:k\in\mathbb{Z}\}$ are such
that
\begin{equation}\label{e:csum}
\sum_{k=-\infty}^\infty |a_k|^r<\infty \quad\text{for some}\quad r<\alpha,\; 0<r\le 1,
\end{equation}
and if~\eqref{eq:cj1} holds then the linear process $\{Z_j\colon j\in\mathbb{Z}\}$ defined by~\eqref{e:map}
satisfies
\[
\frac{1}{b_n}\sum_{j=1}^{n} \bigl(Z_j-Ac_n\bigr) \overset{d}{\longrightarrow}AX(1)\quad
\text{in}\quad\mathbb{R}, \quad \text{where}\quad A=\sum_{k\in\mathbb{Z}}a_k.
\]
For the case $\alpha\in(0,2)$, \cite{avramtaqqu92} show that if $a_k\ge 0$, $k\in\mathbb{Z}$,
satisfy~\eqref{e:csum} and if additional constraints are imposed  for $\alpha\ge 1$ \citep[see][Theorem
2]{avramtaqqu92}, then~\eqref{eq:cj1} implies
\begin{equation}\label{eq:cm1}
\frac{1}{b_n}\sum_{j=1}^{[nt]}(Z_j-Ac_n)\overset{M_1}{\Longrightarrow} AX(t)
\end{equation}
and that the convergence in~\eqref{eq:cm1} is impossible in the $J_1$-topology. We show in Corollary~\ref{c:avt}
that no additional assumptions are needed for  $\alpha\ge 1$.  It is still not known to what extent one can
relax the condition that all $a_k$ have the same sign to get convergence in~\eqref{eq:cm1} with any topology
weaker than $J_1$. With our approach we reduce this problem to continuity properties of addition in a given
topology (see Section~\ref{s:final}). When $\alpha=2$ then~\eqref{eq:cm1} holds with
 any real constants $a_k$ satisfying \eqref{e:csum} \cite[see, e.g.][and the
references therein]{peligradutev06b,moon08}.

We now consider identically distributed, possibly dependent, random variables  $\{\xi_j\colon j\in\mathbb{Z}\}$
with $\E\xi_1^2=\infty$. Note that if \eqref{a:daa} holds with $\alpha\in (0,2]$ then $\E|\xi_1|^\beta<\infty$
for every $\beta\in(0,\alpha)$, thus condition~\eqref{e:csum} ensures that each $Z_j$ in~\eqref{e:map} is a.s.
converging series, since
\[
\E|Z_j|^r\le \sum_{k\in\mathbb{Z}}|a_k|^r\E|\xi_{j-k}|^r=\E|\xi_{1}|^r \sum_{k\in\mathbb{Z}}|a_k|^r<\infty.
\]
Our main result is the following.
\begin{theorem}\label{th:lincon}
Let a linear process $\{Z_j\colon j\in\mathbb{Z}\}$ be defined by~\eqref{e:map}, where $\{\xi_j\colon
j\in\mathbb{Z}\}$ and $\{a_k:k\in\mathbb{Z}\}$  satisfy~\eqref{a:daa}
 and~\eqref{e:csum} with $\alpha\in(0,2]$.  Assume that $\{b_n,c_n\colon n\in\mathbb{N}\}$ are sequences
satisfying the following conditions:
\begin{equation}\label{a:bn}
b_n\to\infty,\quad\frac{c_n}{b_n}\to 0,\quad  \text{and}\quad\limsup_{n\to\infty}n b_n^{-2}E
(\xi_1^2I(|\xi_1|\le b_n))<\infty,
\end{equation}
there exists $s\ge 1$ such that
\begin{equation}\label{a:bi}
\limsup_{n\to\infty}\sup_{k} \E (\max_{1\le l\le
[nT]}|\frac{1}{b_n}\sum_{j=1}^{l}\bigl(\xi_{j-k}I(|\xi_{j-k}|\le b_n)-c_n\bigr)|^s)<\infty
\end{equation}
for all $T>0$, and  there exists a process $X$ such that
\begin{equation}\label{e:ct}
\frac{1}{b_n}\sum_{j=1}^{[nt]}(\xi_{j}-c_n)\Longrightarrow X(t)
\end{equation}
in $\D$ with the topology $J_1$ or $M_1$. If the constants $a_k$ are nonnegative or the process $X$ has
continuous sample paths, then the linear process $\{Z_j\colon j\in\mathbb{Z}\}$ satisfies~\eqref{eq:cm1} with
the same process $X$. Moreover, assumption~\eqref{a:bi} can be omitted, if $\alpha<1$ and $\limsup_{n\to\infty}
nb_n^{-1}|c_n|<\infty$.
\end{theorem}

\begin{remark}
Note that if $b_n\to\infty$ then
\[
\frac{1}{b_n}\E(|\xi_1|I(|\xi_1|\le b_n))=\int_{0}^1\Pr(yb_n<|\xi_1|\le b_n)dy\to0, \quad \text{as}\quad
n\to\infty,
\]
by Lebesgue's dominated convergence theorem, since $\Pr(yb_n<|\xi_1|\le b_n)\to 0$ as $n\to \infty$. Observe
also that if \eqref{e:xitail1} holds with $\alpha\in (0,1)\cup(1,2)$ and the sequences $\{b_n,c_n\colon
n\in\mathbb{N}\}$ are as in \eqref{d:bncn} then 
\[
\frac{1}{b_n}\sum_{j=1}^{[nt]}(\xi_{j}-c_n)\Longrightarrow X(t) \quad \text{if and only if } \quad
\frac{1}{b_n}\sum_{j=1}^{[nt]}(\xi_{j}-c)\Longrightarrow X(t)+\widetilde{c}t,
\]
where $c=0$ for $\alpha<1$,  $c=\E\xi_1$ for $\alpha>1$, and $\widetilde{c}$ is the limit of $nb_n^{-1}(c_n-c)$,
which exists and is finite. The equivalence is also valid when \eqref{a:daa} holds with $\alpha=2$, $b_n$
is such that $nb_n^{-2}\ell(b_n)\to1$, and $c_n$ as above. 
\end{remark}

The proof of Theorem~\ref{th:lincon} is given in Section~\ref{s:pr}. We now comment on  condition~\eqref{a:bi}.
The choice of $c_n=\E(\xi_1I(|\xi_1|\le b_n))$ might allow us to use known  moment maximal inequalities for
partial sums of random variables with mean zero such
as Doob's maximal inequality for martingales or maximal inequalities for strongly mixing sequences.   
In the i.i.d. case Theorem~\ref{th:lincon} implies the following.

\begin{corollary}\label{c:avt}
Let $\{\xi_j\colon j\in\mathbb{Z}\}$ be a sequence of i.i.d. random variables such that condition~\eqref{eq:cj1}
holds, where $X$ is an $\alpha$-stable process with $\alpha\in(0,2)$. If the nonnegative coefficients
$\{a_k:k\in\mathbb{Z}\}$ satisfy~\eqref{e:csum}, then the linear process $\{Z_j\colon j\in\mathbb{Z}\}$
satisfies~\eqref{eq:cm1} with the same process $X$.
\end{corollary}
\begin{proof} We can choose the sequences $\{b_n,c_n\colon
n\in\mathbb{N}\}$ as in \eqref{d:bncn}. Then \eqref{a:bn} holds. For each $n$ and $k\in\mathbb{Z}$, define
\[
\zeta_{n,k,j}:=\xi_{j-k}I(|\xi_{j-k}|\le
  b_n)-\E(\xi_{1}I(|\xi_{1}|\le
b_n)),\quad j\in\mathbb{Z}.
\]
From Doob's maximal inequality for martingales it follows that 
\[
\begin{split}
\E (\max_{1\le l\le [nT]}|\sum_{j=1}^{l}\zeta_{n,k,j}|^2)\le 2\E\bigl(\sum_{j=1}^{[nT]}\zeta_{n,k,j}\bigr)^2
=2[nT]\E(\zeta_{n,k,j})^2\le 2nT\E(\xi_{1}^2I(|\xi_{1}|\le b_n)).
\end{split}
\]
Since $nb_n^{-\alpha}L(b_n)\to 1$ as $n\to\infty$, we have $nb_n^{-2}\E(\xi_1^2I(|\xi_1|\le b_n))\to
\alpha/(2-\alpha)$. Therefore, \eqref{a:bi} holds.
\end{proof}

\begin{remark}
In Corollary~\ref{c:avt}, the convergence in~\eqref{eq:cm1}  can not be strengthened to the $J_1$ topology, by
Theorem~1 of~\cite{avramtaqqu92}. This can also be derived from Theorem~2.4 of \cite{davisresnick85} and
Theorem~3.1 of~\cite{tyran09} \cite[see][Remark 3.3, for more details]{tyran09}.
\end{remark}

Recall that a sequence $\{\xi_j\colon j\in\mathbb{Z}\}$  is said to be \emph{$\rho$-mixing}, if
$\rho(n)\to\infty$, where
\[
\rho(n)=\sup\{|\mathrm{corr}(f,g)|\colon f\in L^2(\F_k), g\in L^2(\F^{n+k}),k\in \mathbb{Z}\}
\]
and $\F_{l}$ ($\F^{l}$) denotes the $\sigma$-algebra generated by $\xi_j$ with indices $j\le l$ ($j\ge l$).

\begin{corollary}\label{c:ms}
Let the sequence $\{\xi_j\colon j\in\mathbb{Z}\}$ be strictly stationary and $\rho$-mixing with $\sum_{i\ge 1}
\rho(2^i)<\infty$. Suppose that condition~\eqref{a:daa} holds with $\alpha\in[1,2]$,  and that the sequence
$\{b_n:n\in\mathbb{N}\}$ satisfies~\eqref{a:bn}, where $c_n=\Pr(\xi_1I(|\xi_1|\le b_n))$, $n\ge 1$. Then
condition~\eqref{a:bi} holds for all $T>0$ with $s=2$.

Moreover, if the coefficients $\{a_k:k\in\mathbb{Z}\}$ satisfy~\eqref{e:csum}, then \eqref{eq:cj1}
implies~\eqref{eq:cm1} with the same process $X$ provided that either the $\{a_k:k\in\mathbb{Z}\}$ are
nonnegative or $X$ has continuous sample paths.
\end{corollary}
\begin{proof}
Note that for each $n$ and $k$ the sequence $\{\zeta_{n,k,j}\colon j\in\mathbb{Z}\}$ defined in the proof of
Corollary~\ref{c:avt} is a stationary sequence of square integrable random variables with mean zero and  is
$\rho$-mixing with at least the same mixing rate as $\{\xi_j\colon j\in\mathbb{Z}\}$. By Theorem~1.1
of~\cite{shao95}, there exists a constant $C>0$, depending only on $\rho(\cdot)$, such that
\[
\E(\max_{1\le l\le [nT]}\bigl(\sum_{j=1}^{l}\zeta_{n,k,j}\bigr)^2)\le C nT\E \zeta_{n,k,1}^2.
\]
Since $\E \zeta_{n,k,1}^2\le \E(\xi_1^2I(|\xi_1|\le b_n))$, condition~\eqref{a:bn} implies~\eqref{a:bi} and the
result follows from Theorem~\ref{th:lincon}.
\end{proof}

In the setting of Corollary~\ref{c:ms}, if we assume that $\rho(1)<1$, then~\eqref{a:daa} with $\alpha=2$
implies that~\eqref{eq:cj1} holds with a sequence $\{b_n:n\in\mathbb{N}\}$ satisfying~\eqref{a:bn} and with $X$
being a standard Brownian motion~\cite[see][]{shao93}. Hence, we recover the result of \cite{moon08}. If
$\alpha\in(0,2)$ then~\eqref{e:xitail1} in general does not imply~\eqref{eq:cj1} as the example of linear
processes shows; see \cite{tyran09} for sufficient conditions when~\eqref{e:xitail1} implies~\eqref{eq:cj1}. In
particular, if in Corollary~\ref{c:ms} we take $\alpha\in[1,2)$, replace~\eqref{a:daa} with~\eqref{e:xitail1},
and choose the sequences $\{b_n,c_n\colon n\in\mathbb{N}\}$ as in~\eqref{d:bncn}, then Theorem 1.1
of~\cite{tyran09} shows that condition~\eqref{eq:cj1} holds with $X$ being an $\alpha$-stable L\'evy process if
and only if for any $\varepsilon>0$ there exist sequences of integers
$r_n=r_n(\varepsilon),l_n=l_n(\varepsilon)\to\infty$ such that
\begin{equation*}\label{eq:asneg0}
r_n=o(n),\quad l_n=o(r_n),\quad n\rho(l_n)=o(r_n), \quad \text{as }n\to \infty,
\end{equation*}
 and
\begin{equation*}\label{eq:asneg}
\lim_{n\to\infty}\Pr(\max_{2\le j\le r_n}|\xi_j|> \varepsilon b_n \bigl||\xi_1|>\varepsilon b_n)=0.
\end{equation*}

\section{Proof of Theorem \ref{th:lincon}}\label{s:pr}

We need the following maximal inequality which follows from Theorem 1 of  \cite{kounias}.
\begin{lemma}\label{l:min} Let $\tau\in(0,1]$.
If $\zeta_1,\ldots,\zeta_N$ are random variables with $E|\zeta_j|^\tau<\infty$ for $j=1,\ldots,N$, then, for any
$\delta>0$,
\[
\Pr(\max_{1\le l\le N}|\zeta_1+\ldots+\zeta_l|> \delta)\le \frac{\sum_{j=1}^N\E(|\zeta_j|^\tau)}{\delta^\tau}.
\]
\end{lemma}

The next lemma will allow us to use Theorem 4.2 of~\cite{billingsley68}. For $m,n\in\mathbb{N}$, define
\[
X_n^{(m)}:=\frac{1}{b_n}\sum_{j=1}^{[nt]}\sum_{|k|\le m}a_k\bigl(\xi_{j-k}-c_n \bigr),\quad
X_n(t):=\frac{1}{b_n}\sum_{j=1}^{[nt]}(Z_j-Ac_n),\quad t\ge 0.
\]

\begin{lemma}\label{l:appr}
Assume \eqref{a:daa}, \eqref{e:csum}, and \eqref{a:bn}. If condition~\eqref{a:bi} holds with $s\ge 1$ and $T>0$
then
\begin{equation}\label{e:prob}
\lim_{m\to\infty}\limsup_{n\to\infty}\Pr(\sup_{0\le t\le T}|X_n(t)-X_n^{(m)}(t)|>\delta)=0\quad \text{for
all}\quad\delta>0.
\end{equation}
If $\alpha<1$ and $\limsup\limits_{n\to\infty}nb_n^{-1}|c_n|<\infty$, then \eqref{e:prob} holds for all $T>0$.
\end{lemma}
\begin{proof}
Define $\xi_{n,j}=\xi_{j}I(|\xi_{j}|\le b_n)-c_n$, $j\in\mathbb{Z},n\in\mathbb{N}$. First note that
\begin{equation*}
\begin{split}
X_n(t)-X_n^{(m)}(t) &=\frac{1}{b_n}\sum_{j=1}^{[nt]}\sum_{|k|> m}a_k
\xi_{n,j-k}+\frac{1}{b_n}\sum_{j=1}^{[nt]}\sum_{|k|> m}a_k \xi_{j-k}I(|\xi_{j-k}|> b_n).
\end{split}
\end{equation*}
Therefore, the probability in \eqref{e:prob} is less than
\begin{equation}\label{e:prob1}
\Pr\bigl(\max_{1\le l\le[nT]}|\frac{1}{b_n}\sum_{j=1}^{l}\sum_{|k|> m}a_k
\xi_{n,j-k}|>\frac{\delta}{2}\bigr)+\Pr\bigl(\max_{1\le l\le[nT]}|\frac{1}{b_n}\sum_{j=1}^{l}\sum_{|k|> m}a_k
\xi_{j-k}I(|\xi_{j-k}|> b_n)|>\frac{\delta}{2}\bigr).
\end{equation}
We now find an upper bound for the first term in~\eqref{e:prob1}. Let $s\ge 1$ be such that condition
\eqref{a:bi} holds. By H\"older's inequality, we have
\[
\bigl(\sum_{|k|> m}|a_k| |\sum_{j=1}^{l}\xi_{n,j-k}|\bigr)^s\le \bigl(\sum_{|k|> m}|a_k|\bigr)^{s-1}\sum_{|k|>
m}|a_k||\sum_{j=1}^{l}\xi_{n,j-k}|^s,
\]
and therefore
\[
\E(\max_{1\le l\le [nT]}|\sum_{j=1}^{l}\sum_{|k|> m}a_k \xi_{n,j-k}|^s)\le \bigl(\sum_{|k|> m}|a_k|\bigr)^s
\sup_{|k|>m} \E (\max_{1\le l\le [nT]}|\sum_{j=1}^{l}\xi_{n,j-k}|^s),
\]
which, by Markov's inequality, leads to
\[
\Pr\bigl(\max_{1\le l\le[nT]}|\frac{1}{b_n}\sum_{j=1}^{l}\sum_{|k|> m}a_k \xi_{n,j-k}|>\frac{\delta}{2}\bigr)\le
\frac{2^s}{\delta^s b_n^s}\bigl(\sum_{|k|> m}|a_k|\bigr)^s\sup_{|k|>m}\E (\max_{1\le l\le
[nT]}|\sum_{j=1}^{l}\xi_{n,j-k}|^s).
\]
From assumption \eqref{a:bi} we conclude that there exists a constant $C_1$  such that, for any $m\ge 1$, we
have
\begin{equation}\label{e:prob2}
\limsup_{n\to\infty}\Pr\bigl(\max_{1\le l\le[nT]}|\frac{1}{b_n}\sum_{j=1}^{l}\sum_{|k|> m}a_k
\xi_{n,j-k}|>\frac{\delta}{2}\bigr)\le C_1 \bigl(\sum_{|k|> m}|a_k|\bigr)^s.
\end{equation}

To estimate the second term in \eqref{e:prob1}, we consider separately the case of $\alpha\in(1,2]$ and
$\alpha\in(0,1]$. Let us note that \eqref{a:daa} and Karamata's theorem imply \citep[see][Lemma, p. 579]{feller71}
\begin{equation*}
\lim_{x\to\infty}\frac{x^{2-\beta}\E(|\xi_1|^\beta I(|\xi_1|>x))}{\E(\xi_1^2 I(|\xi_1|\le x))}=
\frac{2-\alpha}{\alpha-\beta}
\end{equation*}
for all $\beta<\alpha$,  which combined with \eqref{a:bn} gives
\begin{equation}\label{e:karam1}
\limsup_{n\to\infty}n b_n^{-\beta}\E(|\xi_1|^\beta I(|\xi_1|>b_n))<\infty \quad \text{for}\quad \beta<\alpha.
\end{equation}
We  first assume that $\alpha>1$. We have
\[
\E|\sum_{|k|> m}a_k \xi_{j-k}I(|\xi_{j-k}|> b_n)|\le \E(|\xi_{1}|I(|\xi_{1}|> b_n)) \sum_{|k|> m}|a_k|<\infty.
\]
From Lemma~\ref{l:min} it follows that
\[
\Pr\bigl(\max_{1\le l\le[nT]}|\frac{1}{b_n}\sum_{j=1}^{l}\sum_{|k|> m}a_k \xi_{j-k}I(|\xi_{j-k}|>
b_n)|>\frac{\delta}{2}\bigr)\le \frac{2nT}{\delta b_n}\E(|\xi_{1}|I(|\xi_{1}|> b_n)) \sum_{|k|> m}|a_k|.
\]
Applying \eqref{e:karam1} with $\beta=1$, we can find a constant $C_2$ such that
\begin{equation*}\label{e:prob3}
\limsup_{n\to\infty}\Pr\bigl(\max_{1\le l\le[nT]}|\frac{1}{b_n}\sum_{j=1}^{l}\sum_{|k|> m}a_k
\xi_{j-k}I(|\xi_{j-k}|> b_n)|>\frac{\delta}{2}\bigr)\le C_2\sum_{|k|> m}|a_k|,
\end{equation*}
which combined with~\eqref{e:prob1} and \eqref{e:prob2} gives
\[
\limsup_{n\to\infty}\Pr(\sup_{0\le t\le T}|X_n(t)-X_n^{(m)}(t)|>\delta)\le C_1\bigl(\sum_{|k|>
m}|a_k|\bigr)^s+C_2\sum_{|k|> m}|a_k|.
\]
This shows that \eqref{e:prob} holds when $\alpha>1$, since the series $\sum_{k}|a_k|$ converges.

Next, assume that $\alpha\le 1$. Let $r<\alpha$ be as in \eqref{e:csum}.  Since $r\le 1$, we have
\[
\E|\sum_{|k|> m}a_k \xi_{j-k}I(|\xi_{j-k}|> b_n)|^r\le \E(|\xi_{1}|^rI(|\xi_{1}|> b_n)) \sum_{|k|> m}|a_k|^r.
\]
Applying again Lemma~\ref{l:min} and \eqref{e:karam1} with $\beta=r$, we can find a constant $C_3$ such that
\begin{equation*}\label{e:prob4}
\limsup_{n\to\infty}\Pr\bigl(\max_{1\le l\le[nT]}|\frac{1}{b_n}\sum_{j=1}^{l}\sum_{|k|> m}a_k
\xi_{j-k}I(|\xi_{j-k}|> b_n)|>\frac{\delta}{2}\bigr)\le C_3 \sum_{|k|> m}|a_k|^r,
\end{equation*}
which combined with \eqref{e:prob2} completes the proof of~\eqref{e:prob} under the assumption that~\eqref{a:bi}
holds.

To prove the second part of the lemma it suffices to check that~\eqref{e:prob2} remains valid with $s=1$, if
$\alpha<1$ and $\sup_{n}nb_n^{-1}|c_n|<\infty$. We have
\[
\E|\sum_{|k|> m}a_k \xi_{n,j-k}|\le \sum_{|k|> m}|a_k|\E|\xi_{n,j-k}|\le (\E(|\xi_1|I(|\xi_1|\le
b_n))+|c_n|)\sum_{|k|> m}|a_k|.
\]
Therefore, we can apply Lemma~\ref{l:min}, which gives
\[
\Pr\bigl(\max_{1\le l\le[nT]}|\frac{1}{b_n}\sum_{j=1}^{l}\sum_{|k|> m}a_k \xi_{n,j-k}|>\frac{\delta}{2}\bigr)\le
\frac{2nT}{\delta b_n}(\E(|\xi_1|I(|\xi_1|\le b_n)+|c_n|)\sum_{|k|> m}|a_k|.
\]
Since $\alpha<1$, we obtain, by \eqref{a:daa} and  Karamata's theorem \citep[see][Lemma, p. 579]{feller71},
\begin{equation*}
\lim_{n\to\infty}\frac{x\E(|\xi_1| I(|\xi_1|\le x))}{\E(\xi_1^2 I(|\xi_1|\le x))}= \frac{2-\alpha}{1-\alpha}.
\end{equation*}
This together with \eqref{a:bn} shows that $\limsup_{n\to\infty} nb_n^{-1}\E(|\xi_1|I(|\xi_1|\le b_n)<\infty$,
which completes the proof.
\end{proof}

We shall now recall \cite[Chapter 12]{whitt02} that the sequence $\{\psi_n:n\in\mathbb{N}\}$ converges to $\psi$
as $n\to\infty$ in~$\D$ with some Skorohod topology if and only if the restrictions of $\psi_n$ to $[0,T]$
converge to the restriction of $\psi$ to $[0,T]$ in $\DT$ with the same topology for all $T>0$ that are
continuity points of $\psi$. 
Each of the Skorohod topologies in $\DT$ is metrizable with a metric $d_T$  such that $(\DT,d_T)$ is a separable
metric space and weaker than the uniform metric
\begin{equation*}
d_T(\psi_1,\psi_2)\le \sup_{0\le t\le T}|\psi_1(t)-\psi_2(t)| \quad \text{for}\quad \psi_1,\psi_2\in D[0,T].
\end{equation*}

\begin{lemma}\label{p:lch}
Suppose that $\{\xi_j:j\in\mathbb{Z}\}$, $\{b_n,c_n\colon n\in\mathbb{N}\}$, and $X$ are such that
$b_n\to\infty$, $b_n^{-1}c_n\to 0$, and ~\eqref{e:ct} holds
in $\D$ with one of the Skorohod topologies. Then for each $k\in \mathbb{Z}$
\begin{equation}\label{eq:cj1k}
\frac{1}{b_n}\sum_{j=1}^{[nt]}(\xi_{j-k}-c_n)\Longrightarrow X(t),\quad \text{as}\quad n\to\infty,
\end{equation}
in $\D$ with the same topology.
\end{lemma}
\begin{proof}
Let $k\in \mathbb{Z}$.  Define $h_{n}\colon \D\to\D$ by $h_{n}(\psi)(t)=\psi(s_{n}(t))$ for $\psi\in\D$, where
$s_n(t)=\max\{0,t-kn^{-1}\}$. Since $h_{n}(\psi)\to \psi$ in $\D$ with any of the Skorohod topologies,
\eqref{e:ct} implies that
\[
\frac{1}{b_n}\sum_{j=1}^{[ns_n(t)]}(\xi_{j}-c_n)\Longrightarrow X(t),
\]
by Theorem 5.5 of \cite{billingsley68}. If $k<0$, we have
\[
\sup_{t\ge 0}|\frac{1}{b_n}\sum_{j=1}^{[nt]}(\xi_{j-k}-c_n)-\frac{1}{b_n}\sum_{j=1}^{[ns_n(t)]}(\xi_{j}-c_n)|\le
\frac{1}{b_n}|\sum_{j=1}^{1+k}(\xi_{j}-c_n)|\to 0
\]
in probability as $n\to\infty$, since $b_n^{-1}\xi_j\to 0$ in probability for each $j$ and $b_n^{-1}c_n\to 0$.
Similarly, if $k>0$, we have
\[
\sup_{t\ge 0}|\frac{1}{b_n}\sum_{j=1}^{[nt]}(\xi_{j-k}-c_n)-\frac{1}{b_n}\sum_{j=1}^{[ns_n(t)]}(\xi_{j}-c_n)|\le
\frac{1}{b_n}|\sum_{j=1-k}^{0}(\xi_{j}-c_n)|\to 0
\]
in probability as $n\to\infty$. Consequently, the result follows from Slutsky's theorem
\cite[Theorem~4.1]{billingsley68}.
\end{proof}

\begin{lemma}\label{p:cs}
Let $\{\xi_j:j\in\mathbb{Z}\}$  and $\{b_n,c_n\colon n\in\mathbb{N}\}$ be arbitrary sequences. Suppose that
there exists a process $X$  such that for each $k\in\mathbb{Z}$ condition~\eqref{eq:cj1k} holds in $\D$ with the
topology $J_1$ or $M_1$. Then
\begin{equation}\label{a:fapp}
\frac{1}{b_n}\sum_{j=1}^{[nt]}\sum_{|k|\le m}a_k\bigl(\xi_{j-k}-c_n \bigr){\Longrightarrow}\bigl(\sum_{|k|\le m}
a_k\bigr) X(t),\quad \text{as}\quad n\to\infty,
\end{equation}
in $\D$ with the  $M_1$ topology for all constants $a_k\ge 0$, $|k|\le m$, and all $m\ge 0$.
\end{lemma}
\begin{proof}
The constants $a_k$ are nonnegative, thus any two limiting processes $a_{k_1}X$ and $a_{k_2}X$ have jumps of
common sign. Since addition is continuous in the $M_1$ topology for limiting processes with this property
\citep[see e.g.][Theorem 12.7.3]{whitt02}, the result follows.
\end{proof}

\begin{remark}\label{r:cs} Note that addition is continuous in the $J_1$ topology if the limiting processes almost surely have no common
discontinuities. Consequently, if  $X$ has continuous sample paths and the convergence in~\eqref{eq:cj1k} holds
in $\D$ with the $J_1$ topology, then \eqref{a:fapp} holds  with the same topology for any real constants $a_k$,
$|k|\le m$, and all $m\ge 0$.
\end{remark}

\begin{proof}[Proof of Theorem~\ref{th:lincon}]  Lemma~\ref{p:lch} combined with Lemma~~\ref{p:cs} and Remark~\ref{r:cs}
implies that $X_n^{(m)}(t)\overset{M_1}{\Longrightarrow} A_m X(t)$ as $n\to\infty$ for all $m\ge 0$, where
$A_m=\sum_{|k|\le m}a_k$. We also have $A_m X(t)\overset{M_1}{\Longrightarrow}AX(t)$, which is a consequence of
$A_m \to A$ as $m\to \infty$. We can clearly assume that $A\neq 0$ when $X$ has discontinuous sample paths.
Observe that for each $m$ sufficiently large the processes $A_mX$ being a constant multiple of $AX$ have the
same set of continuity points $T_X=\{t>0:\Pr(X(t)= X(t-))=1\}$ with $[0,\infty)\setminus T_X$ being at most
countable. From Lemma~\ref{l:appr} it follows that
\[
\lim_{m\to\infty}\limsup_{n\to\infty}\Pr(d_T(X_n,X_n^{(m)})>\delta)=0\quad \text{for all}\quad \delta>0
\]
for all $T>0$, where $d_T$ is the metric in $\DT$ which induces the $M_1$ topology. Therefore, by Theorem 4.2 of
\cite{billingsley68}, we conclude that 
$X_n(t)\overset{M_1}{\Longrightarrow} A X(t)$ in $\DT$ for all $T\in T_X$, which completes the proof.

\end{proof}

\section{Final remarks}\label{s:final}

Observe that the nonnegativity of the coefficients $\{a_k:k\in\mathbb{Z}\}$ was only used to deduce
\eqref{a:fapp} from \eqref{e:ct}. Thus, we have the following result in one of the Skorohod topologies
$J_1,M_1,J_2,M_2$.


\begin{theorem}\label{th:lincon2}
Let a linear process $\{Z_j\colon j\in\mathbb{Z}\}$ be defined by~\eqref{e:map}, where $\{\xi_j\colon
j\in\mathbb{Z}\}$ and $\{a_k:k\in\mathbb{Z}\}$  satisfy~\eqref{a:daa}
 and~\eqref{e:csum} with $\alpha\in(0,2]$.  Assume that $\{b_n,c_n\colon n\in\mathbb{N}\}$ are sequences
satisfying \eqref{a:bn} and \eqref{a:bi}. 
If there exists a process $X$ such that  for each sufficiently large $m\ge 0$, as $n\to\infty$, \eqref{a:fapp}
holds
in $\D$ with one of the Skorohod topologies, then 
\begin{equation*}
\frac{1}{b_n}\sum_{j=1}^{[nt]}(Z_j-Ac_n)\Longrightarrow AX(t), \quad \text{where}\quad
A=\sum_{k\in\mathbb{Z}}a_k,
\end{equation*}
in $\D$ with the same topology. Moreover,~\eqref{a:bi} can be omitted, if $\alpha<1$ and
$\limsup_{n\to\infty}nb_n^{-1}|c_n|<\infty$.
\end{theorem}

After submission of this paper, the author has learned of a recent result by \cite{basrak10} which gives some
sufficient conditions for \eqref{e:ct} in the $M_1$ topology.

\section*{Acknowledgments} This work was supported by the Natural Sciences and Engineering Research Council
(NSERC, Canada), the Mathematics of Information Technology and Complex Systems (MITACS, Canada), and by Polish
MNiSW grant N N201 0211 33. This research was partially carried out when the author was visiting McGill
University. Helpful questions of the referee,  which lead to significant improvements,  are gratefully
acknowledged.

\end{document}